\documentclass{article}
\usepackage{amsthm,amsfonts, amsbsy, amssymb,amsmath,graphicx}
\usepackage[russian]{babel}

\newtheorem{thm}{Теорема}

\begin{document}
\title{ГЕОДЕЗИЧЕСКИЕ ТКАНИ ГИПЕРПОВЕРХНОСТЕЙ}
\author{@ 2008 г., В.~В. Гольдберг,  В.~В. Лычагин\footnote{\textit{New Jersey Institute of Technology,
USA; Tromso University, Tromso, Norway; email:
vladislav.goldberg@gmail.com, lychagin@yahoo.com}}}
\date{}
\maketitle

\begin{abstract}
В  данной работе мы изучаем геометрические структуры, связанные
с геодезическими тканями гиперповерхностей. Мы показываем, что с каждой
геодезической $(n+2)$-тканью гиперповерхностей на $n$-мерном многообразии естественным образом
связаны единственная проективная структура  и, при условии отмеченного слоения,---
единственная аффинная структура. Проективная структура  выделяется требованием, чтобы слои всех
слоений ткани были вполне геодезическими,  а аффинная структура---дополнительным требованием,
чтобы одна из функций ткани была аффинной.

Эти структуры позволяют определить дифференциальные инварианты геодезических тканей, а также
дать геометрически прозрачные ответы на классические вопросы теории тканей, такие как проблема
линеаризации и теорема Гронвалла.
\end{abstract}

\Large{

\section{Введение} В  данной работе мы изучаем геометрические структуры, связанные
с геодезическими тканями гиперповерхностей. Мы показываем, что с каждой
геодезической $(n+2)$-тканью гиперповерхностей на $n$-мерном многообразии естественным образом
связаны единственная проективная структура  и, при условии отмеченного слоения,---
единственная аффинная структура. Проективная структура  выделяется требованием, чтобы слои всех
слоений ткани были вполне геодезическими,  а аффинная структура---дополнительным требованием,
чтобы одна из функций ткани была аффинной.

Эти структуры позволяют определить дифференциальные инварианты геодезических тканей, а также
дать геометрически прозрачные ответы на классические вопросы теории тканей, такие как проблема
линеаризации и теорема Гронвалла.

Эта работа является непосредственным продолжением работ \cite{AGL04}, \cite{GL09a}, \cite{GL09b} авторов по геодезическим тканям на плоскости.
В \cite{P08} аналогичный вопрос о существовании проективных структур рассмотрен другим методом, который, как
пишет сам автор, не является инвариантным. В  данной работе в отличие от \cite{GL09a} мы используем язык
дифференциальных форм, который позволяет значительно упростить формулы и дать явные выражения
для инвариантов геодезической ткани.

\section{Аффинные связности}
Пусть $M=M^n$---гладкое многообразие размерности $n$, $\nabla$---аффинная
связность без кручения в кокасательном расслоении $T^*M$  и
$d_\nabla$---ковариантный дифференциал :
$$
 d_{\nabla}: \Omega ^{1}(M) \rightarrow \Omega^{1} (M) \otimes \Omega^{1}(M).
 $$
Этот дифференциал может быть представлен в виде
$$
 d_{\nabla }=d\oplus d_{\nabla}^{s},
 $$
где $d$---дифференциал де Рама, а
  $$
d_{\nabla }^{s}: \Omega^{1} (M) \rightarrow S^{2} (\Omega^{1}) (M)
$$
--симметрическая часть дифференциала $d_\nabla$.

\section{Геодезические слоения и проективные структуры}
Необращающаяся в нуль дифференциальная 1-форма $\omega \neq 0$ задает слоение коразмерности один,  если $\omega \wedge d\omega =0$. Слои этого слоения будут вполне
геодезическими в связности  $\nabla$ тогда и только тогда когда (см. \cite{GL09a})
\begin{equation}
d_{\nabla}^{s} \omega = \theta \cdot \omega,\label{geod cond}
\end{equation}%
где $\theta$---некоторая дифференциальная 1-форма.

Мы называем функцию $f$  \emph{вполне геодезической} в связности  $\nabla$, если
ее поверхности уровня $f = \operatorname{const}.$ являются
вполне геодезическими в связности  $\nabla$, и функцию $f$ мы называем
\emph{аффинной} в связности  $\nabla$, если
\begin{equation}
d_{\nabla}^{s}\, df = 0. \label{affine function}
\end{equation}%

Размерность
пространства решений уравнения (\ref{affine function}) будем называть \emph{аффинным рангом} связности  $\nabla$.

Заметим, что аффинный
ранг  связности  $\nabla$ равен размерности многообразия $M$ тогда и только тогда, когда
связность  $\nabla$ плоская.

Две связности  $\nabla$ и $\nabla'$  \emph{проективно эквивалентны}
тогда и только тогда, когда они имеют одни и те же геодезические или (см. \cite{W21})  когда
$$
d_{\nabla }^{s}(\omega)-d_{\nabla '}^{s}(\omega)=\rho \cdot \omega
$$
для некоторой 1-формы $\rho$ и всех 1-форм $\omega$. Более того, соотношение (\ref{geod cond}) показывает, что проективно эквивалентные связности имеют одни и те же вполне геодезические
слоения.

\section{Ткани}
Под \emph{$d$-тканью коразмерности один}
мы понимаем набор $d$ слоений коразмерности один на $M$, если слоения заданы дифференциальными 1-формами $\omega_i, i = 1, \dots , d$, и каждые $n$ из них линейно независимы. Мы обозначим через $<\omega_1, \dots , \omega_d>$
такую $d$-ткань.

Для заданной связности  $\nabla$ мы говорим, что $d$-ткань является
\emph{геодезической}, если слои всех слоений ткани являются вполне геодезическими
в связности  $\nabla$.

Наборы $<\omega_1, \dots , \omega_d>$ и $<s_1 \,\omega_1, \dots , s_d\, \omega_d>$ задают
одну и ту же $d$-ткань, если $s_1 \neq 0, \dots, s_d \neq 0$, где $s_i \in C^\infty (M)$.

 Пусть формы $\omega_1, \dots , \omega_d$ задают $d$-ткань. Выберем $n$ из них, скажем, $\omega_1, \dots , \omega_n$,  за базис. Тогда формы $\omega_{i}, i \geq n+1$, в базисе $\omega_1, \dots , \omega_n$ запишутся в
следующем виде:
\begin{equation}
a_{i1} \omega_1+ \dots + a_{in}\omega_n + \omega_{i} = 0, \label{omega(i)}
\end{equation}%
где координаты $a_{ij}$ не обращаются в ноль, а $a_{i\,n+2}$ мы в дальнейшем для простоты обозначим просто $a_i$. Выбором множителей $s_i, i = 1, \dots , n,$ можно добиться того, что имеет место формула
\begin{equation}
\omega_1+ \dots + \omega_n + \omega_{n+1} = 0. \label{web normalization}
\end{equation}%

В дальнейшем мы будем использовать нормировки (\ref{omega(i)}) и(\ref{web normalization}). Отметим, что при нормировках (\ref{omega(i)}) и (\ref{web normalization}) формы
$\omega_1, \dots , \omega_n, \omega_{n+1}, \dots, \omega_d$ и $s_1 \,\omega_1, \dots ,$
$s_n \,\omega_n, s_{n+1}\, \omega_{n+1}, \dots, s_d \,\omega_d$ задают одну и ту же
$d$-ткань тогда и только тогда, когда $s = s_1 = \dots s_n = s_{n+1} = \dots = s_d$,
а точки $a^{(i)}=[a_{i1}:\dots:a_{in}]$ проективного пространства $\mathbb{RP}^{n-1}$
 являются инвариантами ткани. Мы называем их базисными (см. \cite{AGL04}).

\section{Геодезические ткани}
Из формулы (\ref{geod cond}) вытекает следующий результат:
\begin{thm}
$d$-ткань $<\omega_1, \dots , \omega_d>$ будет геодезической тогда
и только тогда, когда
$$
d_\nabla^s \omega_i = \theta_i \cdot \omega_i, \;\; i = 1, \dots, d,
$$
для некоторых $1$-форм $\theta_i$.
\end{thm}

Выберем базис $\partial_1, \dots , \partial_n$
векторных полей, двойственный кобазису $\omega_1, \dots , \omega_n$: $\omega_i (\partial_j)
= \delta_{ij}$. Тогда
$$
\lbrack\partial_{i},\partial_{j}]=\sum_{k}c_{ij}^{k}\partial_{k}%
$$
для некоторых функций
$c_{ij}^{k}\in C^{\infty} (M)$ и
$$
\nabla_{\partial_i} (\partial_{j})
=\sum_{k}\Gamma_{ji}^{k}\partial_{k}, \;\; 1 \leq i, j \leq n,
$$
где $\Gamma_{ij}^{k}$---символы Кристоффеля второго
рода связности $\nabla$.

Симметричная часть $d_\nabla^s$ принимает вид
\begin{equation}
d_{\nabla}^s (\omega^{k})
=-\sum _{i,j}\Gamma_{ij}^{k}\omega^{j} \cdot \omega^{i}. \label{sympart}
\end{equation}%
Отсюда следует, что
\begin{equation}
\Gamma_{ji}^{k}-\Gamma_{ij}^{k}=c_{ij}^{k}. \label{Gamma-c}
\end{equation}%

Исследуем условия полной геодезичности первых $n+1$ слоений ткани. Для слоений,
определяемых формами  $\omega_1, \dots , \omega_n$, эти условия имеют вид
$$
d_\nabla^s \omega_i = \theta_i \cdot \omega_i, \;\; i = 1, \dots, n,
$$
где
\begin{equation}
\theta_i = \sum_{j=1}^n  \theta_{ij} \omega_j. \label{theta-i}
\end{equation}%
Сравнивая соотношения (\ref{sympart}) и (\ref{theta-i}), получаем
 \begin{equation*}
 \begin{array}{ll}
\Gamma^k_{ik} + \Gamma^k_{ki} + \theta_{ki} = 0, \;\; i, k = 1, \dots, n,\\
\Gamma^k_{ij} + \Gamma^k_{ji} = 0, \;\; \text{если} \;\; i\neq k \;\; \text{и} \;\;j\neq k.
\end{array}
\end{equation*}%
Отсюда и из соотношения(\ref{Gamma-c}) вытекают следующие соотношения между $\Gamma^k_{ij}$ и $c^k_{ij}, \theta_{ki}$:
 \begin{equation}
 \Gamma^k_{ik} = \displaystyle\frac{c^k_{ki} - \theta_{ki}}{2},
\label{Gamma-kik}
\end{equation}%
\begin{equation}
\Gamma^k_{ij} = \displaystyle\frac{c^k_{ji}}{2}, \;\; \text{если} \;\; i\neq k \;\; \text{и} \;\;j\neq k.
\label{Gamma-kij}
\end{equation}%

Обозначим через $\sigma_{ij}$ и $\alpha_{ij}$ симметричную и кососимметричную часть матрицы
$(\theta_{ij})$, т.е.
\begin{equation*}
 \sigma_{ij} = \displaystyle\frac{\theta_{ij} + \theta_{ji}}{2}, \;\;
\alpha_{ij} = \displaystyle\frac{\theta_{ij}-\theta_{ji}}{2},
\end{equation*}
и положим $t_i = \theta_{ii}$.
Тогда условие полной геодезичности $(n+1)$-го слоения полностью определяет
симметрическую часть $\theta_{ij}$,
$$
\sigma_{ij}=\frac{t_i + t_j}{2},
$$
а также дает следующее представление дифференциальной формы $\theta_{n+1}$:
$$
\theta_{n+1}=\sum_{i=1}^n t_i \theta_i.
$$

Условия полной геодезичности $(n+2)$-го слоения позволяют определить кососимметричную
 часть $\alpha_{ij}$:
\begin{equation*}
 \alpha_{ij} = \displaystyle\frac{t_{j} -t_{i}}{2} + s_{ij},
\end{equation*}
где
$$
s_{ij} =s_{ij}^a = \frac{1}{a_i - a_j}\Bigl(a_i  \partial_j - a_{j} \partial_i\Bigr)\log \frac{a_j}{a_i},
 $$
и
\begin{equation*}
 \theta_{n+2} =  \theta_{n+1}
 +  \sum_{i=1}^n \displaystyle\frac{a_{i,i}}{a_i} \omega_i,
\end{equation*}
где $a_{i,i}$--производная от $a_i$ вдоль $\partial_i$.

Окончательно, дифференциальные формы $\theta_{i}, \theta_{n+1}, \theta_{n+2}$ имеют вид
 \begin{equation*}
 \left\{
 \begin{array}{ll}
\theta_{i} = \theta_{n+1} + \displaystyle\sum_{i=1}^n s_{ij} \omega_j,  \;\; i=1, \dots, n,\\
\theta_{n+1} = \theta_{n+1},\\
\theta_{n+2}= \theta_{n+1}+  \displaystyle\sum_{i=1}^n \displaystyle\frac{a_{i,i}}{a_i} \omega_i.
\end{array}
\right.
\end{equation*}%

Отсюда вытекают следующие результаты:
\begin{thm}
Аффинная связность без кручения, для которой $(n+2)$-ткань гиперповерхностей является геодезической,
задается формами $\theta_1, \dots, \theta_n$ вида
\begin{equation}
 \theta_{i} = \theta_{n+1} + \sum_{i=1}^n s_{ij} \omega_j,
\label{forms of aff connection}
\end{equation}%
а соответствующие символы Кристоффеля вычисляются по формулам $(\ref{Gamma-kik})$ и
$(\ref{Gamma-kij})$.
\end{thm}
\begin{thm}\label{pr str}
Каждая
$(n+2)$-ткань гиперповерхностей определяет единственную проективную структуру,
а именно, класс проективно эквивалентных связностей, определяемых формами
$(\ref{forms of aff connection})$.
\end{thm}

Указанную в теореме единственную проективную структуру назовем
\emph{канонической}.

Заметим, что для любой геодезической $d$-ткани, где $d\geq n+2$, канонические проективные
 структуры, определяемые различными $(n+2)$-подтканями, совпадают. Поэтому в дальнейшем мы
   говорим о канонической проективной структуре геодезической $d$-ткани, $d\geq n+2$.

  Ткань с выделенным слоением будем называть \emph{отмеченной}.

Рассмотрим отмеченную  $d$-ткань и предположим, что выделено $(n+1)$-е слоение.
Выберем нормировку (локально) таким образом, чтобы $\omega_{n+1}=df$, а аффинную связность так, чтобы
 форма $\theta_{n+1}\equiv 0$.

Тогда для этой аффинной связности функция $f$ является аффинной. Отметим, что такая функция
$f$ определена с точностью до аффинного калибровочного преобразования $f \rightarrow af+b$.

\begin{thm}\label{aff str}
 Каждая отмеченная
$(n+2)$-ткань гиперповерхностей определяет единственную аффинную связность, для которой ткань является геодезической, а выделенное слоение задается аффинной функцией.
\end{thm}

Указанную в теореме единственную аффинную структуру мы также назовем
\emph{канонической}.

\section{Условия геодезичности $d$-ткани}
Предположим, что мы уже провели нормировки (\ref{web normalization}) и (\ref{omega(i)}). Рассмотрим слоение, задаваемое формой $\omega$, где
$$
\omega = b_1 \omega_1 + \dots + b_n \omega_n.
$$
Это слоение вполне геодезично в связности $\nabla$, если
$$
d_\nabla^s \omega = \theta \cdot \omega,
$$
или
$$
s_{ij}^a = s_{ij}^b.
$$

Отсюда вытекает следующий результат.
\begin{thm}
Обозначим через $(a^k)$ набор базисных инвариантов, где $k = n+2, \dots , d$.
Тогда $d$-ткань гиперповерхностей будет геодезической в том и только том случае, когда
$$
s_{ij}^{(a^k)} = s_{ij}^{(a^l)} \;\; \text{для всех}\;\; k, l = n+2, \dots , d.
$$
\end{thm}

\section{Линеаризуемость тканей}

Известно, что если $\dim M=2$, то необходимым и достаточным условием того, чтобы многообразие $M$ было плоским,
  является обращение в нуль тензорa Лиувилля (см. \cite{Lie83}, \cite{Lio89} или \cite{GL09b}). Если  же $\dim M>2$, то
необходимым и достаточным условием того, чтобы многообразие $M^n$ было плоским, является обращение в нуль
тензора Вейля (см. \cite{VT26}). Отсюда и из результатов раздела 4 вытекает следующая теорема.

\begin{thm} $($\rm{\cite{GL09a}}$)$ $1.$ \emph{Если $\dim M=2$, то $d$-ткань гиперповерхностей,
$d \geq 4$, локально линеаризуема тогда и только тогда, когда она является геодезической, а
 тензор Лиувилля  канонической проективной структуры обращается в нуль.}

$2.$ \emph{Если $\dim M>2$, то $d$-ткань при
$d \geq n+2$ локально линеаризуема тогда и только тогда, когда она являет
ся геодезической, а
  тензор Вейля  канонической  проективной структуры обращается в нуль.}
\end{thm}

\section{Теоремы типа Гронвалла}

Из Теоремы \ref{pr str}
вытекает следующая теорема типа Гронвалла (см. \cite{gr12} и \cite{GL09a} для
$n=2$).

\begin{thm}
Любое отображение геодезической $d$-ткани гиперповерхностей при $d\geq n+2$ на
другую геодезическую $d$-ткань является проективным
преобразованием относительно канонических проективных структур.
\end{thm}

Теорема \ref{aff str}
влечет более сильную теорему типа Гронвалла.

\begin{thm}
Отображение отмеченных геодезических $d$-тканей гиперповерхностей при $d\geq n+2$ является аффинным
относительно канонических аффинных структур.
\end{thm}

\newpage

\begin{center}
\large
GEODESIC WEBS OF HYPERSURFACES
\end{center}

\begin{center}
\bfseries V.\,V.~Goldberg, V.\,V.~Lychagin
\end{center}

\small{
In the present paper we study geometric structures associated with webs of hypersurfaces.
We prove that with any geodesic $(n+2)$-web  on an $n$-dimensional manifold there is naturally associated  a unique projective
structure and, provided that one of web foliations is pointed, there is also associated a unique
affine structure. The projective structure can be chosen by the claim that the leaves of all web foliations
are totally geodesic, and the affine structure by an additional claim that one of web functions
is affine.

These structures allow us to determine differential invariants of geodesic webs and
give  geometrically clear answers to some classical problems of the web theory such as the web
linearization and the Gronwall theorem.}

}


\begin{thebibliography}{9}

\bibitem {AGL04} M.A. Akivis, V.V. Goldberg, V.V. Lychagin,
Selecta Math. \textbf{10}(4), 431--451 (2004).




\bibitem {GL09a} V.~V.~Goldberg,  V.~V.~Lychagin, arXiv: 0810.5392v1, pp. 1--15 (2009) (принято к печати в Acta Appl. Math. (2009)).

\bibitem{GL09b} Goldberg,~V.~V., Lychagin,~V.~V., arXiv: 0812.0125v2, pp. 1--31 (2009) (принято к печати и будет опубликовано в \emph{The Abel Symposium} 2008, Springer (2009)).

\bibitem{gr12} Gronwall, T. H.,
J. de Liouville \textbf{8}, 59--102  (1912).


\bibitem{Lie83}
Lie, S.,   Archiv f\"{u}r Math. og Naturvidenskab \textbf{8} (Kristiania, 1883), 371--458; see also Gesammelte Abhandlungen. Bd. 5 (1924), paper XIV, 362--427.


\bibitem{Lio89} Liouville,~R., Journal de l'\'{E}cole Polytechnique \textbf{59}, 7--76 (1889).



\bibitem{NS94} K.~Nomizu and T.~Sasaki, \textit{Affine Differential Geometry}
(Cambridge Tracts in Mathematics, \textbf{111}. Cambridge
University Press, Cambridge, 1994).




\bibitem{P08} Pirio,~L., arXiv: 0811.1810v1, pp. 1--26 (2008).




\bibitem{VT26} Veblen, O. and Thomas, J. M., Ann. Math. (2) \textbf{27}, no. 3, 279--296, (1926).


\bibitem{W21} Weyl, H., G\"{o}tt. Nachr., 1921, 99--122 (1921).



\end{thebibliography}
\end{document}